# New Vacca-Type Rational Series for Euler's Constant and Its "Alternating" Analog $\ln\frac{4}{\pi}$

## Jonathan Sondow

**1. INTRODUCTION.** In [**11**] we observed that Euler's constant $\gamma$ and $\ln(4/\pi)$ are related by the series

$$\gamma = \sum_{n=1}^{\infty} \left( \frac{1}{n} - \ln\frac{n+1}{n} \right) \tag{1}$$

and

$$\ln\frac{4}{\pi} = \sum_{n=1}^{\infty} (-1)^{n-1} \left( \frac{1}{n} - \ln\frac{n+1}{n} \right), \tag{2}$$

so that $\ln(4/\pi)$ is an "alternating" analog of $\gamma$. In this note we give new rational series for both constants, as well as new proofs of two known series for $\gamma$.

The following expression for Euler's constant is known as *Vacca's series* [**13**]:

$$\gamma = \sum_{n=1}^{\infty} (-1)^n \frac{\lfloor \log_2 n \rfloor}{n}, \tag{3}$$

where $\lfloor x \rfloor$ denotes the greatest integer $\leq x$, and $\log_q n = \frac{\ln n}{\ln q}$ is the logarithm of $n$ to the base $q$. We give an analogous rational series for $\ln(4/\pi)$. (A non-rational analog is given in [**8**, Théorème 2].) We state and (in section 2) prove both formulas simultaneously, using the notation $\gamma^+ = \gamma$ and $\gamma^- = \ln(4/\pi)$. For example, we write

$$\gamma^{\pm} = \sum_{n=1}^{\infty} (\pm 1)^{n-1} \left( \frac{1}{n} - \ln\frac{n+1}{n} \right) \tag{$4_{\pm}$}$$

to mean both (1) and (2).

**Theorem 1.** *If $N_0(n)$ (respectively, $N_1(n)$) is the number of zeros (respectively, ones) in the binary expansion of n, then $\gamma = \gamma^+$ and $\ln(4/\pi) = \gamma^-$ are given by the rational series*

$$\gamma^{\pm} = \sum_{n=2}^{\infty} (-1)^n \frac{N_1\left(\lfloor \frac{n}{2} \rfloor\right) \pm N_0\left(\lfloor \frac{n}{2} \rfloor\right)}{n}. \tag{$5_{\pm}$}$$

Since for $n \geq 2$

$$N_1(\lfloor n/2 \rfloor) + N_0(\lfloor n/2 \rfloor) = \lfloor \log_2 n \rfloor,$$

series $(5_+)$ *for Euler's constant is indeed Vacca's series* (3). Series $(5_\pm)$ begin

$$\gamma = \frac{1}{2} - \frac{1}{3} + \frac{2}{4} - \frac{2}{5} + \frac{2}{6} - \frac{2}{7} + \frac{3}{8} - \frac{3}{9} + \frac{3}{10} - \frac{3}{11} + \frac{3}{12} - \frac{3}{13} + \frac{3}{14} - \frac{3}{15} + \frac{4}{16} - \frac{4}{17} + \cdots,$$

$$\ln\frac{4}{\pi} = \frac{1}{2} - \frac{1}{3} + \frac{2}{6} - \frac{2}{7} - \frac{1}{8} + \frac{1}{9} + \frac{1}{10} - \frac{1}{11} + \frac{1}{12} - \frac{1}{13} + \frac{3}{14} - \frac{3}{15} - \frac{2}{16} + \frac{2}{17} + \frac{2}{22} - \frac{2}{23} + \cdots.$$

**Corollary 1.** *We have the two series*

$$\gamma^\pm = \sum_{n=1}^\infty \frac{N_1(n) \pm N_0(n)}{2n(2n+1)}. \tag{$6_\pm$}$$

*Proof.* In each of the series $(5_\pm)$, group the terms in pairs. ●

We would like to accelerate both series $(6_\pm)$. However, the numerators in series $(6_-)$ for $\ln(4/\pi)$ are irregular (see [**10**, sequence A037861]), and we only succeed in finding accelerated series for $\gamma$.

A. W. Addison used an integral representation of the Riemann zeta function to obtain the series [**1**]

$$\gamma = \frac{1}{2} + \sum_{n=1}^\infty \frac{\lfloor \log_2 2n \rfloor}{2n(2n+1)(2n+2)}, \tag{7}$$

which converges faster than series $(6_+)$ for $\gamma$. In section 3 we prove the following generalization.

**Theorem 2.** *If $q$ is an integer $\geq 2$ and $P_q(x)$ is the polynomial of degree $q-2$*

$$P_q(x) = (qx+1)(qx+2)\cdots(qx+q-1) \sum_{m=1}^{q-1} \frac{m(q-m)}{qx+m}, \tag{8}$$

*then Euler's constant is given by the rational series*

$$\gamma = \frac{1}{2} + \sum_{n=1}^\infty \frac{\lfloor \log_q qn \rfloor P_q(n)}{qn(qn+1)\cdots(qn+q)}. \tag{9}$$

*Moreover, as $q$ increases the series converges more rapidly.*

**Examples.** Since $P_2(x) = 1$, the case $q = 2$ is Addison's series (7). When $q = 3$ we get the faster series

$$\gamma = \frac{1}{2} + \sum_{n=1}^{\infty} \frac{\lfloor \log_3 3n \rfloor (12n + 6)}{3n(3n+1)(3n+2)(3n+3)}.$$

To prove Theorem 2, we employ an averaging technique that S. Krämer used to give a simple proof of Addison's formula by accelerating Vacca's series [**6**, p. 60]. Instead, we accelerate the generalized Vacca series

$$\gamma = \sum_{n=1}^{\infty} \varepsilon(n) \frac{\lfloor \log_q n \rfloor}{n}, \tag{10}$$

where

$$\varepsilon(n) = \begin{cases} q - 1 & \text{if } q \text{ divides } n, \\ -1 & \text{otherwise}. \end{cases}$$

Series (10) was proposed by L. Carlitz [**3**]. It also follows from an integral formula for Euler's constant due to Ramanujan, as B. C. Berndt and D. C. Bowman showed [**2**]. For more on this, as well as other accelerated series for $\gamma$, see [**12**].

In 1897 N. Nielsen discovered a series for $\gamma$ closely related to $(6_+)$ [**7**]. Series (3) was obtained independently by E. Jacobsthal in 1906 [**5**], G. Vacca in 1909 [**13**], H. F. Sandham in 1949 [**9**], and I. Gerst in 1969 [**4**]. Gerst also rediscovered Nielsen's series, and used it to give a new proof of Addison's formula. Papers on Vacca-type series by these and other authors, including J. W. L. Glaisher (1910), G. H. Hardy (1912), J. C. Kluyver (1924), V. Brun (1938), S. Selberg (1940, 1967), R. W. Gosper, Jr. (1972), M. Koecher (1989), and F. L. Bauer (1990), are discussed in [**2**, p. 23] and [**6**, sections 5.2.2 and 10.1].

**2. PROOF OF THEOREM 1.** We first indicate the proofs of three lemmas.

**Lemma 1.** *Let*

$$A_n = \frac{1}{n} - \ln \frac{n+1}{n}. \tag{11}$$

*Then*

$$A_n = \frac{1}{2n(2n+1)} + A_{2n} + A_{2n+1} \tag{12}$$

*and*

$$0 < A_n < \frac{1}{n} - \frac{1}{n+1}. \tag{13}$$

*Proof.* Relation (12) follows directly from (11). To prove inequalities (13), use the integral formula

$$A_n = \int_n^{n+1} \left( \frac{1}{n} - \frac{1}{x} \right) dx. \qquad \bullet$$

**Lemma 2.** *If*
$$\Delta^{\pm}(n) = N_1(n) \pm N_0(n),$$
*then for* $n \geq 2$
$$\Delta^{\pm}(\lfloor n/2 \rfloor) + (\pm 1)^{n-1} = \Delta^{\pm}(n).$$

*Proof.* Verify the formulas for even $n$ and odd $n > 1$. ●

**Lemma 3.** *Let*
$$S_k^{\pm} = \sum_{n=1}^{k} (\pm 1)^{n-1} A_n$$

*be the k-th partial sum of series* $(4_{\pm})$ *for* $\gamma^{\pm}$. *Then*

$$S_{2^k-1}^{\pm} = \sum_{n=1}^{2^{k-1}-1} \frac{\Delta^{\pm}(n)}{2n(2n+1)} + R_k^{\pm}, \qquad (14_{\pm})$$

*where the remainder term is*

$$R_k^{\pm} = \sum_{n=2^{k-1}}^{2^k-1} \Delta^{\pm}(n) A_n. \qquad (15_{\pm})$$

*Proof.* We induct on $k$. When $k = 1$ we have $S_1^{\pm} = A_1 = R_1^{\pm}$, as required (since the sum in $(14_{\pm})$ is empty). Now write

$$S_{2^{k+1}-1}^{\pm} = S_{2^k-1}^{\pm} + \sum_{n=2^k}^{2^{k+1}-1} (\pm 1)^{n-1} A_n,$$

and invoke the inductive hypotheses $(14_{\pm})$ and $(15_{\pm})$ after substituting (12) into $(15_{\pm})$. Use the identity

$$\sum_{n=2^{k-1}}^{2^k-1} \Delta^{\pm}(n)(A_{2n} + A_{2n+1}) = \sum_{n=2^k}^{2^{k+1}-1} \Delta^{\pm}(\lfloor n/2 \rfloor) A_n$$

and Lemma 2 to complete the inductive step. ●

*Proof of Theorem 1.* The proof of Corollary 1 from Theorem 1 can be reversed, so it suffices to prove that

$$\gamma^{\pm} = \sum_{n=1}^{\infty} \frac{\Delta^{\pm}(n)}{2n(2n+1)}. \qquad (16_{\pm})$$

Since by its definition $S_k^{\pm} \to \gamma^{\pm}$ as $k \to \infty$, we can deduce $(16_{\pm})$ using $(14_{\pm})$ if $R_k^{\pm} \to 0$ as $k \to \infty$. But by $(15_{\pm})$, (13), and telescoping,

$$\left|R_k^{\pm}\right| \leq \sum_{n=2^{k-1}}^{2^k-1}\left|\Delta^{\pm}(n)\right|A_n < \sum_{n=2^{k-1}}^{2^k-1} k\left(\frac{1}{n}-\frac{1}{n+1}\right) = \frac{k}{2^k} \to 0.$$   ●

**3. PROOF OF THEOREM 2.** If we group the terms of series (10) in $q$-tuples with positive first member, we obtain

$$\gamma = \sum_{n=1}^{\infty}\left\lfloor \log_q qn \right\rfloor\left(\frac{q-1}{qn}-\frac{1}{qn+1}-\frac{1}{qn+2}-\cdots-\frac{1}{qn+q-1}\right). \qquad (17)$$

If instead we group the terms in $q$-tuples with positive last member, we deduce that

$$\gamma = 1 + \sum_{n=1}^{\infty}\left\lfloor \log_q qn \right\rfloor\left(-\frac{1}{qn+1}-\frac{1}{qn+2}-\cdots-\frac{1}{qn+q-1}+\frac{q-1}{qn+q}\right), \qquad (18)$$

where the first 1 is the result of summing the geometric series

$$\frac{q-1}{q}+\frac{q-1}{q^2}+\frac{q-1}{q^3}+\cdots = 1.$$

Now average (17) and (18) to get

$$\gamma = \frac{1}{2} + \frac{1}{2}\sum_{n=1}^{\infty}\left\lfloor \log_q qn \right\rfloor\left(\frac{q-1}{qn}-\frac{2}{qn+1}-\frac{2}{qn+2}-\cdots-\frac{2}{qn+q-1}+\frac{q-1}{qn+q}\right).$$

The expression in parentheses is equal to

$$\sum_{m=1}^{q-1}\left(\frac{1}{qn}-\frac{2}{qn+m}+\frac{1}{qn+q}\right) = \frac{1}{qn(qn+q)}\sum_{m=1}^{q-1}\left(2m-q+\frac{2m(q-m)}{qn+m}\right),$$

and $\sum_{m=1}^{q-1}(2m-q) = 0$, so

$$\gamma = \frac{1}{2}+\sum_{n=1}^{\infty}\frac{\left\lfloor \log_q qn \right\rfloor}{qn(qn+q)}\sum_{m=1}^{q-1}\frac{m(q-m)}{qn+m}.$$

In view of (8), this gives the desired formula (9).

In order to show that as $q$ increases series (9) converges faster, we fix $q$ and compute that as $n \to \infty$ the $n$-th term of the series is asymptotic to the product

$$\frac{1-q^{-2}}{\ln q}\cdot\frac{\ln n}{6n^3}.$$

The first factor is a decreasing function of $q$, and the proof of Theorem 2 is complete.

**ACKNOWLEDGMENTS.** I am grateful to Stefan Krämer and Wadim Zudilin for useful comments on an early version of this note, and to Tanguy Rivoal for sending me a draft of [**8**].

*209 West 97th Street, New York, NY 10025*
*jsondow@alumni.princeton.edu*